\documentclass[11pt, a4]{article}
\usepackage[english]{babel}
\usepackage[utf8x]{inputenc}
\usepackage{hyphenat}
\usepackage{amsmath}
\usepackage{mathrsfs}
\usepackage{color}
\usepackage{url}
\usepackage{graphicx}
\usepackage{multicol}
\usepackage[bottom]{footmisc}
\usepackage{tabularx, booktabs, rotating}
\usepackage{multirow}
\usepackage{subfigure}
\usepackage{natbib}
\usepackage{changepage}
\usepackage{changepage}
\usepackage{algorithm}
\usepackage{algpseudocode}
\usepackage[normalem]{ulem}
\usepackage{mathtools}
\usepackage{enumerate}
\usepackage[shortlabels]{enumitem}
\usepackage{adjustbox}
\usepackage{hyperref}

\algnewcommand{\To}{{\normalfont\bfseries to }}
\algnewcommand{\Continue}{{\normalfont\bfseries continue}}

\usepackage{pgfplots}
\pgfplotsset{width=10cm,compat=1.9}

\date{}

\title{Solution of a Practical Vehicle Routing Problem for Monitoring Water Distribution Networks}
\author{Reza Atefi$^{(1)}$,  Manuel Iori$^{(2)}$, Majid Salari$^{(1)}$ and Dario Vezzali$^{(3)}$\\ \\
 $(1)$ Department of Industrial Engineering, Ferdowsi University of Mashhad, \\ P.O. Box 91779-48951, Mashhad, Iran \\
 \texttt{atefi@mail.um.ac.ir, msalari@um.ac.ir} \\
 $(2)$ Department of Sciences and Methods for Engineering,\\ University of Modena and Reggio Emilia, \\ Via Amendola 2, 42122, Reggio Emilia, Italy \\
 \texttt{manuel.iori@unimore.it} \\
 $(3)$ ``Marco Biagi" Foundation, University of Modena and Reggio Emilia, \\ Largo Marco Biagi 10, 41121, Modena, Italy \\
 \texttt{dario.vezzali@unimore.it} \\
}

\begin{document}
\maketitle

\begin{abstract}
\noindent In this work, we introduce a generalization of the well-known Vehicle Routing Problem for a specific application in the monitoring of a Water Distribution Network (WDN). In this problem, multiple technicians must visit a sequence of nodes in the WDN and perform a series of tests to check the quality of water. Some special nodes (i.e., wells) require technicians to first collect a key from a key center. The key must then be returned to the same key center after the test has been performed, thus introducing precedence constraints and multiple visits in the routes. To solve the problem, three mathematical models and an Iterated Local Search have been implemented. The efficiency of the proposed methods is demonstrated by means of extensive computational tests on randomly created instances, as well as on instances derived from a real-world case study.
\end{abstract}


\section{Introduction} \label{sec:introduction}
Water contamination is related to the presence of one or more chemical compounds or pathogens to the extent that they become dangerous to the consumer and might lead to diseases \citep{NNM2018}. The risk of accidental contamination of drinking water is a well-known issue, and, recently, concerns regarding the deliberate contamination of urban water networks have called for additional safeguards. 

In general, any threat to urban water networks directly affects the users in the community \citep{R2015}. Indeed, according to a report recently released by the \textit{World Health Organization}, contaminated drinking water is estimated to cause 485,000 diarrhoeal deaths each year \citep{WHO2019}. The safety of water distribution networks has always been an important issue for the communities. However, many distribution systems in cities around the world face the threat of accidental or intentional contamination during the transportation from treatment plants to consumers due to reverse flows (i.e., the return of contaminated water flows from facilities), old infrastructures, insufficient use of disinfectants, and so forth. Consequently, water contamination in distribution networks is considered as the most diffused cause behind the spread of water-borne diseases \citep{MKCS2002}.

In recent years, several studies have been conducted to identify the main sources of water pollution and improve the quality of water thanks to innovative treatment methods and plants, but still an accidental event, such as a large-scale contamination or a destructive attack to the transmission system, can significantly affect both the economy and the society. In 2014, for example, 300,000 consumers in West Virginia were affected by the accidental contamination of their drinking water distribution system caused by 4-Methylcyclohexanemethanol \citep{SA2006}. During the same year, as reported by \cite{MDV2017}, a spill of benzene from a chemical plant in China accidentally reached the water distribution network. More recently, 27,000 Norwegian consumers were exposed to water contaminated with Clostridium \citep{deWPWK2019}. 

Supply, treatment, transmission and distribution of drinking water in urban distribution networks require substantial expenses; therefore, not only water in urban distribution networks is considered an essential resource, but also an economic commodity. The results of a study conducted by the \textit{World Bank} show that nearly 15\% of treated water is wasted annually in developed countries. This amount arises to a range of 35-60\% for developing countries \citep{ZWZQHZ2016}. Timely control of \textit{Water Distribution Networks} (WDNs) is thus of fundamental importance, both from an economical and public health point of view.

In this paper, a new variant of the well-known \textit{Vehicle Routing Problem} (VRP) in the context of WDNs is proposed. In this problem, a set of technicians must visit a set of nodes, including wells, reservoirs and treatment plants, within a network to evaluate the water quality. When visiting a well, the technicians need a key to open the well and perform the required tests. Since the technicians do not have the key, they have to visit a specified node at which the key is located, called \textit{key center} in the following, to acquire it. As a result, they need to visit this node before reaching the well. After the tests have been performed, they have to take the key back to its original key center before returning to the depot where they started their route. In addition to that, it is imposed that all nodes are visited and that the duration of any route performed by a technician does not exceed a maximum traveling time. The aim of the problem is to minimize the sum of the traveled times.  

The problem originates from a real-world application that we encountered in Mashhad (Iran), where 5 technicians daily inspect a WDN comprising 3,124 households/shops, 293 reservoirs/tanks, 356 wells and 14 treatment plants. To solve the problem, we propose three \textit{Mixed Integer Linear Programming} (MILP) models, and an \textit{Iterated Local Search} (ILS) algorithm. While the models managed to solve small-size instances with up to 20 nodes, the ILS efficiently tackled cases with up to 200 nodes, allowing us to produce good-quality solutions for randomly created instances, as well as for realistic instances derived from the case study, in short computing times.

The remainder of the paper is organized as follows. In Section \ref{sec:literature_review}, the relevant literature is revised. The problem is formally described in Section \ref{sec:problem_description}. Sections \ref{sec:mathematical_models} and \ref{sec:ILS} present the mathematical models and the ILS algorithm, respectively. Computational results are described in Section \ref{sec:computational_results}, and final conclusions and future research directions are discussed in Section \ref{sec:conclusions}.

\section{Literature Review} \label{sec:literature_review}
The VRP is an iconic class of problems in operations research, with applications in the fields of transportation, distribution, logistics and services. We refer the interested reader to \citet{TV2014} for an extensive overview, and to \citet{MS2022} for a recent survey.
The problem we face generalizes the VRP by considering precedence constraints and multiple visits. In this section, we only revise routing problems involving these two features, with a particular focus on real-world applications.

In the context of the {\textit{Traveling Salesman Problem} (TSP)}, precedence constraints were first addressed in the seminal work by \citet{BFP1995}, and, since then, have been widely investigated.
In \citet{MKCS2002}, the authors proposed a formulation for the TSP with precedence constraints using a two-commodity network flow model and developed a genetic algorithm based on a topological sorting of customers. In \citet{SSB2005}, novel formulations for the asymmetric TSP and the precedence constrained asymmetric TSP were proposed. To tighten the formulations, the authors proposed and tested valid inequalities. \citet{SVDVW2018} presented a new model for the time-dependent capacitated profitable tour problem, a generalization of the TSP with time windows and precedence constraints, and developed a tailored labeling algorithm. \citet{SED2020} describe the precedence constrained generalized TSP, in which customers are partitioned into groups and exactly one visit per group must be performed. They presented a novel branching technique and compared several bounding methods.

Precedence constraints have also been widely studied for problems involving multiple vehicles. \citet{R2015} developed a genetic algorithm based on a topological sorting of customers to solve the VRP with precedence constraints. The algorithm includes a route repair method to generate feasible offspring.
A VRP variant with time windows, synchronization and precedence constraints was introduced by \citet{HLP2016}. The authors focused on an attended home health care application, and proposed some exact and heuristic solution methods, including a novel MILP formulation, a greedy heuristic, and three metaheuristics.

Precedence constraints naturally arise in the context of \textit{Pickup-and-Delivery Problems} (PDP), where each demand must be first collected at an origin node before being delivered at a destination node. We refer the reader to \citet{BCI2014} and \citet{DS2014} for detailed surveys on PDPs for goods transportation and PDPs for people transportation, respectively, and to \citet{KLT2020} for a recent survey on simultaneous PDPs.
Recently, \citet{ACC2020} studied a multi-PDP with time windows. They defined a 2-index formulation, an asymmetric representatives formulation, and a 3-index formulation improved by preprocessing and valid inequalities. The problem was solved exactly using a branch-and-cut algorithm. Dedicated branch-and-cut algorithms were also developed by \citet{HLR2021}, to solve the single-vehicle two-echelon one-commodity PDP, and by \citet{WS2021}, to solve a PDP with split loads and transshipments.
The problem addressed in the latter work includes multiple visits to the same node. This is common when split deliveries are allowed, or multiple pickup and delivery operations can be performed at a single node. These generalizations were considered by \citet{BI2017}, where non-elementary formulations were proposed for a single-vehicle PDP and then extended to the cases of split deliveries, intermediate drop-offs, and multiple vehicles.

Overall, we may find many routing problems that are inspired by real-world applications and involve precedence constraints and multiple visits. \citet{SPS2004} studied an application of a PDP with time windows and precedence constraints arising in the transportation of live animals. In this case, the precedence constraints are given by veterinary rules, imposing that the livestock holdings are visited in a predefined sequence to avoid the spread of potential diseases. The authors proposed a tight formulation of the problem based on a Dantzig-Wolfe decomposition.
\citet{QLLH2013} presented an application in the context of military operations, that was modeled as a generalized VRP with synchronization and precedence constraints. The peculiarity of the problem is due to the nature of the attack, which may require aircraft synchronization, multiple attacks to the same target, and precedence constraints among different targets. The problem was solved by a MILP model.

\citet{FMM2017} addressed a particular PDP with time windows originating from the oil industry. The aim of the problem is to determine the routing and scheduling of vessels that collect crude oil from offshore platforms and transport it to terminals on the coast. The authors proposed a MILP model, solving it by means of two different branch-and-cut algorithms.
Another valuable example of routing and scheduling in the context of large-scale disaster relief operations was examined by \citet{SBMH2019}. The authors solved a PDP arising from a case study in the city of Tehran (Iran). They  proposed an integrated logistic system to evacuate people from areas affected by natural or man-made disasters. The problem was formulated as a MILP model, and a memetic algorithm was developed to solve large-scale instances.
\citet{PAM2020} studied a particular \textit{Workforce Scheduling and Routing Problem} (WSRP), called multiperiod WSRP with dependent tasks, in which the requested services consist of tasks to be executed along one or more days by teams of workers having different skills. Each customer can be visited more than once in a day, as long as precedence constraints are not violated. A MILP model, a constructive algorithm and an ant colony metaheuristic were proposed.

Recently, an interesting variant of the \textit{Team Orienteering Problem} (TOP), named multi-visit TOP with precedence constraints, was investigated by \citet{HMZ2020}. In this problem, a set of tasks has to be accomplished in a predetermined order by possibly different vehicles. To solve the problem, the authors proposed a compact MILP formulation and a kernel search heuristic.

For what concerns WDNs, the literature mainly contains works on the location of sensors (see, e.g., \citealt{RG2014}). The VRP has been applied in many areas, but, to the best of our knowledge, not yet to the inspection of WDNs. In this paper, we fill this lack in the literature and propose exact and heuristic solution methods for a real-world VRP on a WDN.

\section{Problem Description} \label{sec:problem_description}
The WDN is an essential infrastructure that consists of many elements, including reservoirs, wells, pipes and treatment plants.

An effective way to constantly monitor a WDN is by means of water quality sensors, which can be positioned all over the network. In cities where these sensor systems have not been installed, technicians are required to regularly visit nodes of the WDN and perform tests. The nodes to be visited, called for simplicity \textit{demand nodes} in the following, are divided into two types:
\begin{enumerate}
    \item \textit{Type I}: households, shops, reservoirs, tanks and treatment plants. For this kind of nodes, the technicians can directly go on site and perform the required tests. Reservoirs, tanks and treatment plants are characterized by larger service times than households and shops, due to the larger amount of tests that have to be performed;
    \item \textit{Type II}: wells. For these nodes, the technicians need a key to access the well and perform the tests. So, they have to visit first a specified key center, and take the key. Once all tests have been completed at the well, the key needs to be returned to its original key center, thus imposing a second visit.
\end{enumerate}

\begin{figure}
    \centering
    \includegraphics[width=0.85\textwidth]{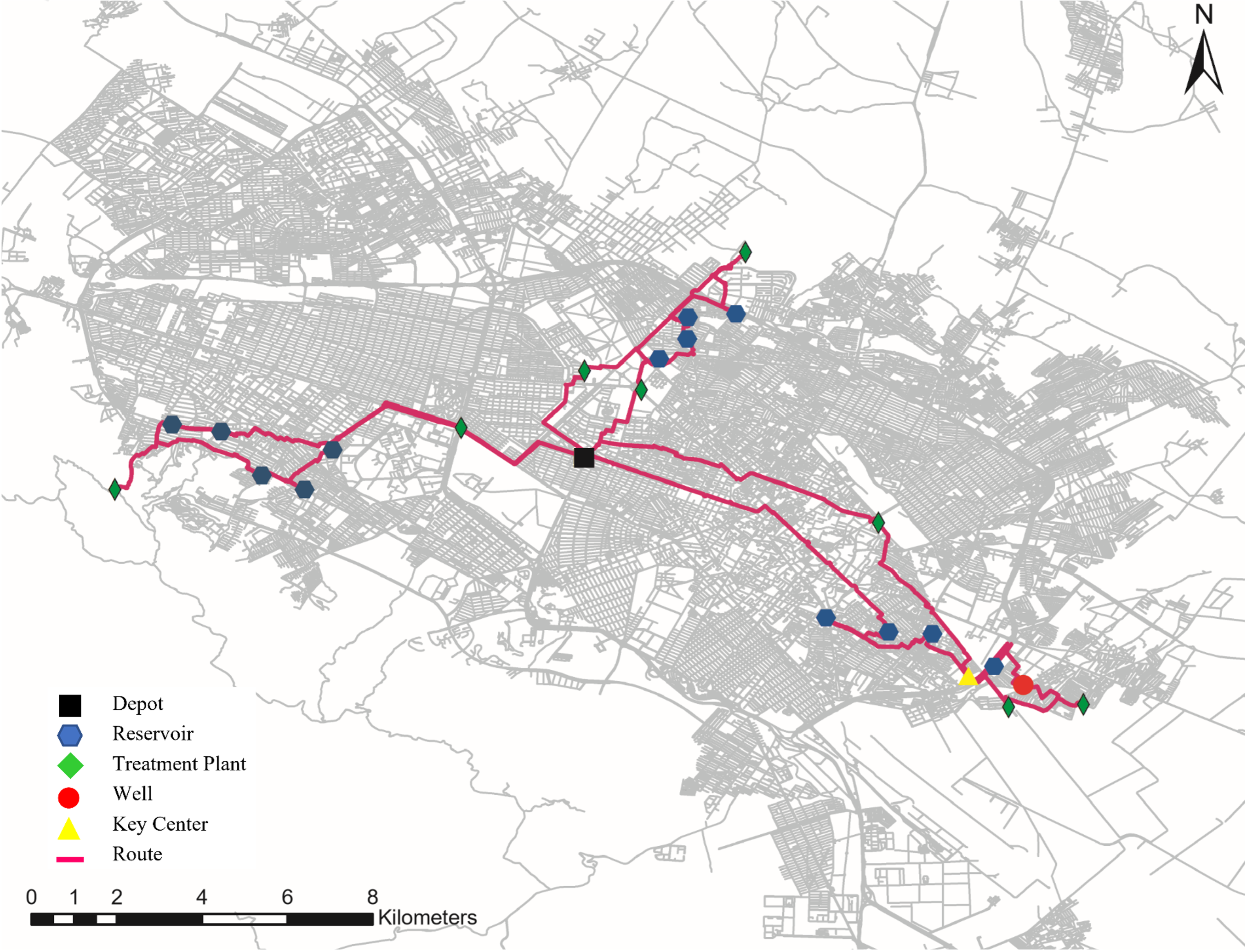}
    \caption{\small An illustrative example of a VRPWDN solution in Mashhad (Iran)}
    \label{fig:example1}
\end{figure}

A simple illustrative example derived from the real-world application we are facing is depicted in Figure~\ref{fig:example1}. It comprises three routes starting and ending at the depot. Two of them (top and left part of the figure) visit just reservoirs and treatment plants, so demand nodes of type I. The third (right part of the figure) also visits a well, and is thus forced to pass twice by the corresponding key center.

Formally, we are given a directed graph $G = (V, A)$, where the node set is $V = \{0,1,\dots,n,n+1\}$ and is partitioned as $V = V_1 \cup V_2 \cup V_3 \cup \{0,n+1\}$. Nodes 0 and $n+1$ represent, respectively, the beginning and end of all routes, and in our application coincide with a unique central depot. Sets $V_1$ and $V_2$ are associated with, respectively, the demand nodes of types I and II. Set $V_3$ comprises nodes associated with all key centers. With each node $i \in V_2$, we associate a predecessor $p_i \in V_3$ and a successor node $d_i \in V_3$. In our application, $p_i$ and $d_i$ correspond to a unique key center, so they have the same geographical location, but the models and algorithms that we propose below can also solve the case in which they correspond to different locations.

Each demand node has to be visited exactly once, while each vehicle visits a particular key center at most once for picking up all the keys, and then another single time for delivering all the keys that were previously collected. This implies that, in case a center has the keys for multiple demand nodes and these nodes are visited by a unique vehicle, then such keys must be collected all together in a unique visit (to $p_i$), and then later delivered all together in another visit (to $d_i$). Note that it is not compulsory to visit $p_i$ immediately before $i$. In other words, the vehicle can collect the key for $i$ but then visit other nodes before reaching $i$. The same holds for $d_i$, which is not required to be visited immediately after $i$.

The graph is complete, and with each arc $(i,j) \in A$ we associate a traveling time $c_{ij}$. A service time $v_i$ is associated with each node $i \in V$. We suppose that triangle inequality holds for all our instances (i.e., $c_{ij} \leq c_{ik} + v_k + c_{kj}$ for all $i, j, k \in V$). We are also given a set $K$ of homogeneous vehicles based at the central depot. Each vehicle performs a single route.

A route starts and ends at the depot. Its duration is given by the sum of the service and traveling times of the nodes and arcs covered by the vehicle, and it should not exceed a maximum duration $L$. Whenever a route visits a node $i$ of type II, then it should also visit $p_i$ and $d_i$. The aim of the \textit{Vehicle Routing Problem for Water Distribution Networks} (VRPWDN) is to visit all demand nodes, while satisfying all constraints and minimizing the sum of the route durations.

The VRPWDN is NP-hard in the strong sense, because it generalizes the well-known VRP. In the next sections, we attempt its solution through mathematical models and heuristic algorithms.

\section{Mathematical Models} \label{sec:mathematical_models}
In this section, we investigate three mathematical models that describe the VRPWDN and are derived from the literature.
The first model is based on a time representation of the problem and is inspired by the formulation proposed by \citet{DMR2014} for the VRP with time windows.
The second is a flow-based model that builds upon the formulation presented by \citet{K2011} and later used by, among others, \citet{KAKD2012}, \citet{NS2014}, and \citet{ASV2015}.
The third is a node-based model that we derive from the classical Miller, Tucker and Zemlin formulation (see, e.g., \citealt{BG2014}).

\subsection{Time-based Model} \label{subsec:time_based_model}
Let $y_{ik}$ be a binary variable taking value 1 if node $i$ is visited by vehicle $k$ and 0 otherwise, $x_{ijk}$ be another binary variable taking value 1 if arc $(i,j)$ is covered by vehicle $k$ and 0 otherwise, and $t_{ik}$ be a continuous variable corresponding to the time at which vehicle $k$ arrives at node $i$.
The time-based model for the VRPWDN can be formulated as follows:
{\small
    \begin{align}
        \label{eq:VRPWDNtb_obj} \text{({VRPWDN}\textsubscript{tb})} \quad & \min~\sum_{k \in K}\sum_{i \in V \setminus \{n+1\}}\sum_{j \in V \setminus \{0\}} (c_{ij} + v_i) x_{ijk}
    \end{align}
}
subject to
{\small
\allowdisplaybreaks
\begin{align}
    \label{eq:VRPWDNtb_cstr2} & \sum_{k \in K}\sum_{j \in V \setminus\{0\}} x_{ijk} = 1 & i \in V_1 \cup V_2 \\
    \label{eq:VRPWDNtb_cstr3} & \sum_{j \in V \setminus\{0\}} x_{0jk} = 1 & k \in K \\
    \label{eq:VRPWDNtb_cstr6} & y_{ik} = \sum_{j \in V \setminus\{0\}} x_{ijk} = \sum_{j \in V \setminus \{n+1\}} x_{jik} & i \in V, k \in K \\
    \label{eq:VRPWDNtb_add_cstr} & \sum_{i \in V \setminus\{n+1\}} x_{i,n+1,k} = 1 & k \in K \\
    \label{eq:VRPWDNtb_cstr7} & \sum_{k \in K} t_{0k} = 0 & \\
    \label{eq:VRPWDNtb_cstr8} & 0 \leq t_{ik} \leq L y_{ik} & i \in V, k \in K \\
    \label{eq:VRPWDNtb_cstr9} & t_{jk} \geq t_{ik} + v_i + c_{ij}  - M_{ij} (1 - x_{ijk}) & \hspace{-2cm} i \in V \setminus \{n+1\}, j \in V \setminus \{0\}, k \in K \\
    \label{eq:VRPWDNtb_cstr11} & y_{p_ik} + y_{d_ik} \geq 2 y_{ik} & i \in V_2, k \in K \\
    \label{eq:VRPWDNtb_cstr12-13b} & t_{p_ik} + v_{p_i} + c_{p_ii} - M_{i}' (1 - y_{ik}) \leq t_{ik} \leq t_{d_ik} - (c_{id_i} + v_i) y_{ik} & i \in V_2, k \in K \\
    \label{eq:VRPWDNtb_dom1} & x_{ijk} \in \{0,1\} & i,j \in V, k \in K \\
    \label{eq:VRPWDNtb_dom2} & y_{ik} \in \{0,1\} & i \in V, k \in K
\end{align}
}
Objective function~\eqref{eq:VRPWDNtb_obj} is to minimize the total duration of the routes. Constraints~\eqref{eq:VRPWDNtb_cstr2} impose that each node $i \in V_1 \cup V_2$ has exactly one outgoing arc. Each vehicle starts its route from the depot and such condition is imposed by means of constraints~\eqref{eq:VRPWDNtb_cstr3}. Constraints~\eqref{eq:VRPWDNtb_cstr6} and \eqref{eq:VRPWDNtb_add_cstr} ensure that each node $i$ has exactly one incoming and one outgoing arc and that each vehicle $k$ end its route at the depot. Constraints~\eqref{eq:VRPWDNtb_cstr7} impose that all routes start at time $0$. Constraints~\eqref{eq:VRPWDNtb_cstr8} impose that arrival times are non-negative and limit the duration of each route to be at most $L$. The time at which vehicle $k$ arrives at node $j$ is modeled by means of constraints~\eqref{eq:VRPWDNtb_cstr9}, in which we set $M_{ij} = L + v_i + c_{ij} - c_{j,n+1}$. Constraints~\eqref{eq:VRPWDNtb_cstr11} impose that if vehicle $k$ visits node $i$, then it also visits nodes $p_i$ and $d_i$. Since $p_i$ may contain keys not only for $i$ but for other nodes, vehicle $k$ may visit $p_i$ but not $i$, and the same holds for $d_i$. For this reason, the equation cannot be an equality. Constraints~\eqref{eq:VRPWDNtb_cstr12-13b}, in which we set $M_{i}' = L + v_{p_i} + c_{p_i, i}- c_{i, n+1}$, guarantee the respect of precedence constraints, by forcing time dependency between visits to $p_i$, $i$ and $d_i$. Note that if $y_{ik}$ is equal to 0, constraints~\eqref{eq:VRPWDNtb_cstr12-13b} become redundant with respect to constraints~\eqref{eq:VRPWDNtb_cstr8}. Constraints~\eqref{eq:VRPWDNtb_dom1} and \eqref{eq:VRPWDNtb_dom2} define the domain of the {$x_{ijk}$ and $y_{ik}$} variables.

Furthermore, the aforementioned model can be enhanced with the addition of the following valid inequalities
{\small
\begin{align}
    \label{eq:VRPWDNtb_vi1} & (c_{0p_i} + v_{p_i} + c_{p_ii}) y_{ik} \leq t_{ik} & i \in V_2, k \in K \\
    \label{eq:VRPWDNtb_vi2} & (c_{0p_i} + v_{p_i} + c_{p_ii} + v_i + c_{id_i}) y_{ik} \leq t_{d_ik} & i \in V_2, k \in K \\
    \label{eq:VRPWDNtb_vi3} & t_{jk} \geq (c_{0i} + v_i + c_{ij}) x_{ijk} & i \in V \setminus \{n+1\}, j \in V \setminus \{0\}, k \in K
\end{align}
}
which strengthen the values taken by the arrival time variables.

\subsection{Flow-based Model} \label{subsec:flow_based_model}
Let $f_{ijk}$ be a variable representing the ``load'' of vehicle $k$ when traveling along arc $(i,j) \in A$. The load represents the number of nodes visited by vehicle $k$ before it travels along arc $(i,j)$.
We can model the VRPWDN as follows:
{\small
\begin{align}
    \label{eq:VRPWDNfb_obj} \text{(VRPWDN\textsubscript{fb})} \quad & \min~\sum_{k \in K}\sum_{i \in V \setminus \{n+1\}}\sum_{j \in V \setminus \{0\}} (c_{ij} + v_i) x_{ijk}
\end{align}
}
subject to \eqref{eq:VRPWDNtb_cstr2}, \eqref{eq:VRPWDNtb_cstr3}, \eqref{eq:VRPWDNtb_dom1} and
{\small
\allowdisplaybreaks
\begin{align}
    \label{eq:VRPWDNfb_cstr22} & \sum_{j \in V \setminus \{0\}} x_{ijk} = \sum_{j \in V \setminus \{n+1\}} x_{jik} & i \in V, k \in K \\
    \label{eq:VRPWDNfb_cstr24} & \sum_{i \in V \setminus \{n+1\}}\sum_{j \in V \setminus \{0\}} (c_{ij} + v_i) x_{ijk} \leq L & k \in K \\
    \label{eq:VRPWDNfb_cstr25} & \sum_{j \in V \setminus \{0\}} f_{0jk} = 0 & k \in K \\
    \label{eq:VRPWDNfb_cstr26} & \sum_{i \in V \setminus \{n+1\}} f_{i,n+1,k} = \sum_{i \in V \setminus \{n+1\}}\sum_{j \in V \setminus \{0\}} x_{ijk} - 1 & k \in K \\
    \label{eq:VRPWDNfb_cstr27} & \sum_{j \in V \setminus \{0\}} f_{ijk} \geq \sum_{j \in V \setminus \{n+1\}} (f_{jik} + x_{jik}) & i \in V \setminus \{0,n+1\}, k \in K \\
    \label{eq:VRPWDNfb_cstr29} & \sum_{j \in V \setminus \{0\}} (f_{p_ijk} - f_{ijk} + x_{ijk}) \leq (n - 1) (1 - \sum_{j \in V \setminus \{0\}} x_{ijk}) & i \in V_2, k \in K \\
    \label{eq:VRPWDNfb_cstr30} & \sum_{j \in V \setminus \{0\}} f_{p_ijk} \geq \sum_{j \in V \setminus \{0\}} x_{ijk} & i \in V_2, k \in K \\
    \label{eq:VRPWDNfb_cstr31} & \sum_{j \in V \setminus \{0\}} (f_{d_ijk} - f_{ijk}) \geq \sum_{j \in V \setminus \{0\}} x_{ijk} & i \in V_2, k \in K \\
    \label{eq:VRPWDNfb_cstr28} & 0 \leq f_{ijk} \leq (n - 1) x_{ijk} & i,j \in V, k \in K
\end{align}
}
As in the previous model, objective function~\eqref{eq:VRPWDNfb_obj} minimizes the total route duration. Constraints~\eqref{eq:VRPWDNfb_cstr22} correspond to the previous constraints~\eqref{eq:VRPWDNtb_cstr6} except for the $y_{ik}$ term. The maximum duration of each route is bounded by means of constraints~\eqref{eq:VRPWDNfb_cstr24}. Constraints~\eqref{eq:VRPWDNfb_cstr25} and \eqref{eq:VRPWDNfb_cstr26} impose the load on vehicle $k$ when leaving 0 and entering $n+1$, respectively.
Constraints~\eqref{eq:VRPWDNfb_cstr27} impose the load conservation at node $i$. Constraints~\eqref{eq:VRPWDNfb_cstr29}--\eqref{eq:VRPWDNfb_cstr31} guarantee the respect of precedence constraints. Constraints~\eqref{eq:VRPWDNfb_cstr28} impose lower and upper bounds on the $f_{ijk}$ variables.

The above model can be improved by the addition of the following constraints:
{\small
\begin{align}
    \label{eq:VRPWDNfb_cstr32} & \sum_{j \in V \setminus \{n+1\}} x_{jp_ik} + \sum_{j \in V \setminus \{0\}} x_{d_ijk} \geq 2 \sum_{j \in V \setminus \{n+1\}} x_{jik} & i \in V_2, k \in K\\
    \label{eq:VRPWDNfb_cstr33} & \sum_{j \in V \setminus \{0\}} (f_{d_ijk} - f_{p_ijk}) \geq \sum_{l \in V \setminus \{n+1\}:p_i = p_l} \sum_{j \in V \setminus \{0\}} x_{ljk} & i \in V_2, k \in K \\
    \label{eq:VRPWDNfb_cstr34} & \sum_{l \in V \setminus \{n+1\}}(f_{lik} - f_{ljk}) + n x_{ijk} + (n - 2) x_{jik} \leq (n - 1) & i,j \in V \setminus \{0,n+1\}, k \in K
\end{align}
}
Constraints~\eqref{eq:VRPWDNfb_cstr32} are equivalent to \eqref{eq:VRPWDNtb_cstr11}. Constraints~\eqref{eq:VRPWDNfb_cstr33} enforce an additional relation between the flows leaving $p_i$ and $d_i$. Constraints~\eqref{eq:VRPWDNfb_cstr34} are derived from the lifted constraints proposed by \citet{DL1991}.

\subsection{Node-based Model} \label{subsec:node_based_model}
Let $u_{ik}$ be a variable representing the load on vehicle $k$ after leaving node $i$. With respect to the previous model, this implies setting $u_{ik}=\sum_{j {\in V}} f_{ijk}$. We can model the VRPWDN as follows:
{\small
    \begin{align}
        \label{eq:VRPWDNnb_obj} \text{(VRPWDN\textsubscript{nb})} \quad & \min~\sum_{k \in K}\sum_{i \in V \setminus \{n+1\}}\sum_{j \in V \setminus \{0\}} (c_{ij} + v_i) x_{ijk}
    \end{align}
}
subject to \eqref{eq:VRPWDNtb_cstr2}, \eqref{eq:VRPWDNtb_cstr3}, \eqref{eq:VRPWDNtb_dom1}, \eqref{eq:VRPWDNfb_cstr22}, \eqref{eq:VRPWDNfb_cstr24} and
{\small
\allowdisplaybreaks
\begin{align}
    \label{eq:VRPWDNnb_cstr41} & u_{0k} = 0 & k \in K \\
    \label{eq:VRPWDNnb_cstr42} & u_{n+1,k} = \sum_{i \in V \setminus \{n+1\}}\sum_{j \in V \setminus \{0\}} x_{ijk} & k \in K \\
    \label{eq:VRPWDNnb_cstr43} & u_{ik} - u_{jk} + n x_{ijk} \leq (n - 1) & i \in V \setminus\{n+1\}, j \in V \setminus\{0\}, k \in K \\
    \label{eq:VRPWDNnb_cstr45} & u_{p_ik} - u_{ik} + \sum_{j \in V \setminus \{0\}} x_{ijk} \leq n (1 - \sum_{j \in V \setminus \{0\}} x_{ijk}) & i \in V_2, k \in K \\
    \label{eq:VRPWDNnb_cstr46} & u_{p_ik} \geq \sum_{j \in \{0\}} x_{ijk} & i \in V_2, k \in K \\
    \label{eq:VRPWDNnb_cstr47} & u_{d_ik} - u_{ik} \geq \sum_{j \in V \setminus \{0\}} x_{ijk} & i \in V_2, k \in K \\
    \label{eq:VRPWDNnb_cstr44} & 0 \leq u_{ik} \leq n \sum_{j \in V \setminus \{0\}} x_{ijk} & i \in V, k \in K
\end{align}
}
For each vehicle $k$, constraints~\eqref{eq:VRPWDNnb_cstr41} set the load after leaving node $0$, while constraints~\eqref{eq:VRPWDNnb_cstr42} define the load when arriving at node $n+1$. Constraints~\eqref{eq:VRPWDNnb_cstr43} impose the load conservation when traveling from node $i$ to node $j$. Constraints~\eqref{eq:VRPWDNnb_cstr45}--\eqref{eq:VRPWDNnb_cstr47} guarantee the respect of precedence constraints. Constraints~\eqref{eq:VRPWDNnb_cstr44} impose both the non-negativity of the $u_{ik}$ variables and their relation with the $x_{ijk}$variables.

The model can be improved by the addition of \eqref{eq:VRPWDNfb_cstr32} and of
{\small
\begin{align}
    \label{eq:VRPWDNnb_cstr49} & u_{d_ik} - u_{p_ik} \geq \sum_{l \in V \setminus \{n+1\} : p_i = p_l}\sum_{j \in V \setminus \{0\}} x_{ljk} & i \in V_2, k \in K \\
    \label{eq:VRPWDNnb_cstr50} & u_{ik} - u_{jk} + n x_{ijk} + (n - 2) x_{jik} \leq (n - 1) & i,j \in V \setminus \{0,n+1\}, k \in K
\end{align}
}
which correspond to the above \eqref{eq:VRPWDNfb_cstr33} and \eqref{eq:VRPWDNfb_cstr34}, respectively.

\section{Iterated Local Search} \label{sec:ILS}
We developed an ILS algorithm with the purpose of finding good-quality VRPWDN solutions in short computing times. The choice of this metaheuristic is motivated by its simplicity and effectiveness, in addition to the wide applicability it has found on related VRPs (see, e.g., \citealt{VSBVO2009}, \citealt{SSO2015}, \citealt{HLP2016} and \citealt{ASCR2018}).
On the other hand, the need for short computing times is justified by the number of visits usually scheduled in a day in our real-world application, and by the fact that candidate locations might change at the beginning or in the course of a day. Two examples which typically cause a re-scheduling of visits can be a new warning for potential water contamination coming from a household or shop or, when visiting a well, the unfortunate event that the well's door is broken and it is not possible to open it.

Following the general framework proposed by \citet{LMS2019}, the ILS starts from an initial solution and then improves it by iteratively invoking local search and perturbation procedures. The pseudo-code of the proposed ILS is provided in Algorithm~\ref{ils:main}. First, we generate an initial solution $x_0$ by means of a heuristic algorithm (line 1), and then we improve it with a local search procedure (line 2). The current solution, $x$, is stored as the incumbent, $x^*$, and inserted in the set of best known solutions obtained during the search, called $BKSet$ (lines 3 and 4). Next, we execute two phases, one after the other.

In the first phase, by applying a perturbation on $x$ followed by a call to the local search (lines 6--8), the algorithm tries to escape from local optima. The perturbation is randomly selected between two tailored procedures. Let $z(x)$ and $l(x)$ be the cost of $x$ and the maximum duration of a route in $x$, respectively. In case $x$ has better cost than $x^*$, or same cost but lower maximum duration, then we use it to update $x^*$. In such a case, we also insert $x$ in $BKSet$. This set contains the $\beta$ different solutions found during the search and having the smallest $z(x)$ costs, breaking ties by smallest $l(x)$ value.
If, instead, $x$ does not improve $x^*$, then we set $x \gets x^*$ as starting solution to be shaken at the next iteration. This loop is repeated until no improvement is found for $max_{iter}$ iterations.

With the aim of further improving the solution obtained, at line 16 we enter the second ILS phase, in which a new series of improving attempts is performed. The idea is to intensify the search around the solutions contained in $BKSet$. For each such solution, we perform once more a loop of shaking and local search procedures, which is repeated until the same termination condition used above is met. Should one of these attempts manage to improve the incumbent solution, this time only in terms of costs, then the search restarts from the beginning of the first phase.

In the following, we provide the details of the main elements of the algorithm.

\begin{algorithm}[htb]
    \caption{\small Iterated Local Search (ILS)}\label{ils:main}
    \algnewcommand{\EndIIf}{\unskip\ \algorithmicend\ \algorithmicif}
    \begin{algorithmic}[1]
        \State $x_0 \gets \texttt{Initialization}()$ \Comment{Generate an initial solution}
        \State $x \gets \texttt{LocalSearch}(x_0)$
        \State $x^* \gets x$
        \State $BKSet \gets \{x^*\}$ \Comment{$BKSet$: set of best known solutions}        
        \Repeat \Comment{Phase 1}
            \State $\texttt{Shake}() \gets \texttt{Rand}\{S_1, S_2\}$ \Comment{Randomly select a shaking procedures}
            \State $x \gets \texttt{Shake}(x)$
            \State $x \gets \texttt{LocalSearch}(x)$
            \State $\texttt{Insert}(x, BKSet)$
            \If{$z(x) < z(x^*)$ OR $(z(x) = z(x^*)$ AND $l(x) < l(x^*))$}
                \State $x^* \gets x$
            \Else
                \State $x \gets x^*$
            \EndIf
        \Until{no improvement is found for $max_{iter}$ iterations}
        \For{$j \gets 1$ \To $|BKSet|$} \Comment{Phase 2}
            \State $x \gets BKSet_j$ \Comment{Select the $j^{th}$ solution $\in BKSet$}
            \Repeat
                \State $\texttt{Shake}() \gets \texttt{Rand}\{S_1, S_2\}$
                \State $x \gets \texttt{Shake}(x)$
                \State $x \gets \texttt{LocalSearch}(x)$
                \If{$z(x) < z(x^*)$}
                    \State $x^* \gets x$
                    \State $\texttt{Insert}(x^*, BKSet)$
                    \State $\texttt{Go to line 5}$
                \EndIf
            \Until{no improvement is found for $max_{iter}$ iterations}
        \EndFor
        \State \Return $x^*$
    \end{algorithmic}
\end{algorithm}

\subsection{Initialization Procedure} \label{subsec:init_procedure}
Algorithm~\ref{ils:init_procedure} gives the \texttt{Initialization} procedure that is used to generate an initial solution. At the beginning, $|K|$ routes are built in parallel by randomly selecting a first node $i \in V_1 \cup V_2$ per route. In case $i$ belongs to $V_2$, then the predecessor and the successor of $i$ (i.e., $p_i$ and $d_i$) are also inserted into the route.
In the next $|V_1 \cup V_2| - |K|$ iterations, a new node is randomly selected and inserted into an existing route. In these iterations, both the node and, in case $i \in V_2$, its predecessor and successor are inserted in the route in the positions that lead to the minimum extra mileage cost. Note that the insertion of node $i$ or tuple $(p_i, i, d_i)$ into an existing route is led by procedure \texttt{CheapestInsertion}, which evaluates among the $|K|$ routes the best candidate for the expansion. At line 22, the algorithm checks whether the solution is feasible. If not, then the whole procedure is repeated from scratch.

\begin{algorithm}[htb]
    \caption{\small Initialization Procedure}\label{ils:init_procedure}
    \begin{algorithmic}[1]
        \State $\mathcal{S},\mathcal{V} \gets \emptyset$
        \For{$k \gets 1$ \To $|K|$} \Comment{Initialization of $|K|$ routes in parallel}
            \State $i \gets \texttt{Rand}\{1, ..., |V_1 \cup V_2|\}$
            \State $\mathcal{V} \gets \mathcal{V} \cup \{i\}$ \Comment{Add $i$ to the set of visited nodes}
            \If{$i \in V_2$}
                \State $r_k \gets (0, p_i, i, d_i, n+1)$
                \State $\texttt{Insert}(r_k, \mathcal{S})$
            \Else
                \State $r_k \gets (0, i, n+1)$
                \State $\texttt{Insert}(r_k, \mathcal{S})$
            \EndIf
        \EndFor
        \For{$j \gets 1$ \To $|V_1 \cup V_2| - |K|$} \Comment{Expansion of existing routes}
            \State $i \gets \texttt{Rand}\{\{1, ..., |V_1 \cup V_2|\} \setminus \mathcal{V}\}$
            \State $\mathcal{V} \gets \mathcal{V} \cup \{i\}$
            \If{$i \in V_2$}
                \State $\texttt{CheapestInsertion}((p_i, i, d_i), r_k \in \mathcal{S})$
            \Else
                \State $\texttt{CheapestInsertion}(i, r_k \in \mathcal{S})$
            \EndIf
        \EndFor
        \If{$\texttt{Feasible}(\mathcal{S}) = 1$}
            \State \Continue
        \Else
            \State \texttt{Go to line 1}
        \EndIf
        \State \Return $\mathcal{S}$
    \end{algorithmic}
\end{algorithm}

\subsection{Local Search} \label{subsec:LS}
The \texttt{LocalSearch} procedure invokes, one after the other, the following neighborhood searches:
\begin{enumerate}[start = 1,label = {LS\arabic*}]
    \item \textit{Swap intra-route}: swap two sequences with up to three consecutive nodes in the same route. Potential nodes belonging to $V_3$ are extracted from the two sequences and reinserted after the swap following the minimum extra mileage cost and respecting the precedence constraints;
    \item \textit{Swap inter-route}: swap two sequences with up to three consecutive nodes from different routes, taking care of nodes belonging to $V_3$;
    \item \textit{Relocate intra-route}: remove a sequence with up to three consecutive nodes and reinsert it in a different position within the same route, taking care of nodes belonging to $V_3$;
    \item \textit{Relocate inter-route}: remove a sequence with up to three consecutive nodes and reinsert it in a different route, taking care of nodes belonging to $V_3$;
    \item \textit{3-opt}: in a preliminary step, select a route and remove potential nodes belonging to $V_3$. Following this step, apply the standard 3-opt algorithm to the remaining nodes. After each iteration of the 3-opt algorithm update the solution by reinserting the previously extracted nodes belonging to $V_3$. 
\end{enumerate}

Procedures from LS1 to LS4 have all complexity $O(n^2)$, whereas LS5 has complexity $O(n^3)$. To limit the computational effort required by LS5, a random logic search is added. In particular, a candidate route $k$ is selected randomly and potential nodes belonging to $V_3$ are removed as follows. For each node $i$ in the route, the saving $s_i$ that could be obtained by removing $i$ and directly connecting the predecessor and successor nodes of $i$ in the route is computed.
Then, the probability of removing $i$ is set to $p_i = {s_i}/{\sum_{{j}}{s_j}}$. By means of the roulette wheel mechanism, three non-adjacent nodes are selected for removal, and then the resulting route is optimized by a 3-opt algorithm. A threshold of $\gamma$ iterations is set to limit the number of attempts.

The calls to LS1--LS5 are repeated as long as an improvement is found. Procedure \texttt{LocalSearch} hence returns a solution which is a local optimum with respect to all five neighborhoods.

\subsection{Shaking Procedure} \label{subsec:shaking}
To perturb a solution, we randomly select, with same probability, one of the two following procedures. 
\begin{enumerate}[start = 1,label = {S\arabic*}]
    \item \textit{Shaking 1:} randomly select a route $k$ and execute a random iteration of the 3-opt algorithm to update the order of visits. If the cost of the current solution is not worse than $\alpha z(x^*)$, with $\alpha$ being an input parameter, randomly select a second route $k'$ and perform another 3-opt iteration. The procedure is iterated as long as the cost of the perturbed solution is not worse than $\alpha z(x^*)$;
    \item \textit{Shaking 2:} compute the cost saving obtained by removing any node from the solution, similarly to what is done in LS5. Then use the roulette wheel mechanism to select a node $i \in V \setminus \{0,n+1\}$, and remove $i$ from its route. The removal procedure is iterated until at least $\alpha$ percent of all nodes have been removed. If the selected node belongs to $V_2$, then its saving is computed as the average cost saving obtained by removing $i$, $p_i$, and $d_i$. At the end of this step, the algorithm invokes the \texttt{Initialization} procedure to rebuild a feasible solution.
\end{enumerate}

\section{Computational Results} \label{sec:computational_results}
In this section, we present the results of extensive computational tests performed with the aim of assessing the performance of the proposed methods. 
The mathematical models and the ILS were coded in C++ using Microsoft Visual Studio 2010. The computational tests were executed on a PC equipped with an Intel Core i7 CPU processor @ 2.70 GHz and 6 GB of RAM, using CPLEX 12.3 as MILP solver.
In Section~\ref{subsec:test_instances}, we describe the sets of randomly-created instances that we used for our tests. The comparison among the mathematical models is reported in Section~\ref{subsec:mathematical_models_comparison}, while the behavior of the ILS is analyzed in Sections~\ref{subsec:ILS_parameter_tuning} and \ref{subsec:ILS_evaluation}. In Section~\ref{subsec:realistic_instances}, we report the results of additional computational experiments performed on a set of realistic instances derived from the case study.

\subsection{Randomly-created Instances} \label{subsec:test_instances}
We created several random instances with the aim of assessing the performance of the algorithms under different situations. In detail, we created two sets of instances, each comprising different subsets having homogeneous values of $|V_1 \cup V_2|$, $(|V_2|,|V_3|)$ and $|K|$, and composed by three random instances per subset. We obtained the following sets:
\begin{itemize}
\item \emph{Small-size}: 
18 instances with $|V_1 \cup V_2|$=10, 
$(|V_2|,|V_3|) \in \{(1,1)$, $(2,1)$, $(2,2)\}$, and $|K| \in \{1,2\}$; 	
24 instances with $|V_1 \cup V_2|$=15, 
$(|V_2|,|V_3|) \in \{(3,2)$, $(3,3)$, $(4,2)$, $(4,3)\}$, and $|K| \in \{2, 3\}$; 
24 instances with $|V_1 \cup V_2|$=20, 
$(|V_2|,|V_3|) \in \{(2,2)$, $(3,2)$, $(3,3)$, $(5,3)\}$, and $|K| \in \{2, 3\}$; 
\item \emph{Medium- and large-size}: 
24 instances with $|V_1 \cup V_2|$=50, 
$(|V_2|,|V_3|) \in \{(5,5)$, $(8,8)$, $(10,5)$, $(10,8)\}$, and $|K| \in \{5, 8\}$; 	
24 instances with $|V_1 \cup V_2|$=100, 
$(|V_2|,|V_3|) \in \{(5,5)$, $(10,5)$, $(10,10)$, $(15,10)\}$, and $|K| \in \{10, 15\}$; 
24 instances with $|V_1 \cup V_2|$=200, 
$(|V_2|,|V_3|) \in \{(10,10)$, $(20,10)$, $(20,20)$, $(30,20)\}$, and $|K| \in \{15, 20\}$.
\end{itemize}

For each instance, the coordinates of the nodes are integer values randomly selected between 0 and 100. The distances between the nodes are computed as the Euclidean ones, rounded to the second closest digit. The maximum duration is set to
$L = 1.5 ({\sum_{i \in V_1 \cup V_2} \overline{c}_i + |K| \sum_{i \in V_3 \cup \{0\}} \overline{c}_i})/{|K|}$, where $\overline{c}_i$ is the average travel time of the arcs leaving $i$, computed as $\overline{c}_i = \sum_{j \in V \setminus\{i\}} c_{ij}/(|V|-1)$ for each node $i \in V \setminus\{n+1\}$. The service time $v_i$ for each node $i \in V_1 \cup V_2 \cup V_3$ is set to a random integer value between 20 and 40.

In the following, a subset of instances is identified by the tuple $(|V_1 \cup V_2|, |V_2|, |V_3|,$ $|K|)$, while a single instance is identified by $(|V_1 \cup V_2|, |V_2|, |V_3|, |K|, u)$, where $u$ is a numerical index going from 1 to 3.

To favor future research on the problem, the randomly-created instances have been made publicly available at \url{https://github.com/DarioVezzali/VRPWDN}.

\subsection{Comparison among the Mathematical Models} \label{subsec:mathematical_models_comparison}
In this section, the performance of the three mathematical models from Section~\ref{sec:mathematical_models} is investigated. A time limit of 3,600 CPU seconds was imposed on each execution. The aggregated results that we obtained are reported in Table~\ref{tab:math_models_comparison}. 
Each line reports average/total values for a group of three instances having the same numbers of vertices and vehicles. 
For each group, columns ``$z_{lb}$'' and ``$z_{ub}$'' give the average lower and upper bound values, respectively, column ``$\%$gap'' gives the average percentage gap and column ``t(s)'' the average run time. An entry ``tlim'' indicates that the time limit was reached for all the three instances in the group. Column ``opt'' gives the total number of instances solved to proven optimality.

\begin{sidewaystable}
    \caption{\small Comparison of mathematical models (three inst. per line). Best average lower and upper bound values in \textbf{boldface}} \label{tab:math_models_comparison}
    \centering
    \scriptsize
    \setlength{\tabcolsep}{2.5pt}
    \adjustbox{max width=\columnwidth}{
        \begin{tabular}{c c c c r r r r r r r r r r r r r r r}
            \toprule
            & & &  & \multicolumn{5}{c}{time-based} & \multicolumn{5}{c}{flow-based} & \multicolumn{5}{c}{node-based} \\
            \cmidrule(lr){5-9} \cmidrule(lr){10-14} \cmidrule(lr){15-19}
            $|V_1 \cup V_2|$ & $|V_2|$ & $|V_3|$ & $|K|$ & $z_{lb}$ & $z_{ub}$ & $\%$gap &    t(s) & opt & $z_{lb}$ & $z_{ub}$ & $\%$gap &    t(s) & opt & $z_{lb}$ & $z_{ub}$ & $\%$gap &    t(s) & opt \\
            \cmidrule(lr){1-4} \cmidrule(lr){5-9} \cmidrule(lr){10-14} \cmidrule(lr){15-19}
            10 & 1 & 1 & 1 &  \textbf{706.39} &  \textbf{706.39} &    0.00 &    0.87 &   3 &  \textbf{706.39} &  \textbf{706.39} &    0.00 &    1.13 &   3 &  \textbf{706.39} &  \textbf{706.39} &    0.00 &    1.42 &   3 \\
            10 & 1 & 1 & 2 &  \textbf{730.91} &  \textbf{730.91} &    0.00 &    7.45 &   3 &  \textbf{730.91} &  \textbf{730.91} &    0.00 &    3.26 &   3 &  \textbf{730.91} &  \textbf{730.91} &    0.00 &   12.20 &   3 \\
            10 & 2 & 1 & 1 &  \textbf{743.25} &  \textbf{743.25} &    0.00 &    6.54 &   3 &  \textbf{743.25} &  \textbf{743.25} &    0.00 &    3.13 &   3 &  \textbf{743.25} &  \textbf{743.25} &    0.00 &    0.66 &   3 \\
            10 & 2 & 1 & 2 &  \textbf{793.75} &  \textbf{793.75} &    0.00 &   23.12 &   3 &  \textbf{793.75} &  \textbf{793.75} &    0.00 &    9.45 &   3 &  \textbf{793.75} &  \textbf{793.75} &    0.00 &   10.96 &   3 \\
            10 & 2 & 2 & 1 &  \textbf{803.96} &  \textbf{803.96} &    0.00 &  177.59 &   3 &  \textbf{803.96} &  \textbf{803.96} &    0.00 &   16.73 &   3 &  \textbf{803.96} &  \textbf{803.96} &    0.00 &   44.80 &   3 \\
            10 & 2 & 2 & 2 &  \textbf{833.78} &  \textbf{833.78} &    0.00 &  110.68 &   3 &  \textbf{833.78} &  \textbf{833.78} &    0.00 &   75.01 &   3 &  \textbf{833.78} &  \textbf{833.78} &    0.00 &  436.67 &   3 \\
            \cmidrule(lr){1-4} \cmidrule(lr){5-9} \cmidrule(lr){10-14} \cmidrule(lr){15-19}
            \multicolumn{4}{l}{sum/avg (10)} & \textbf{768.67} & \textbf{768.67} & 0.00 & 54.38 & 18 & \textbf{768.67} & \textbf{768.67} & 0.00 & 18.12 & 18 & \textbf{768.67} & \textbf{768.67} & 0.00 & 84.45 & 18 \\
            \cmidrule(lr){1-4} \cmidrule(lr){5-9} \cmidrule(lr){10-14} \cmidrule(lr){15-19}
            15 & 3 & 2 & 2 & \textbf{1036.61} & \textbf{1036.61} &    0.00 & 1811.64 &   3 & \textbf{1036.61} & \textbf{1036.61} &    0.00 &  509.71 &   3 & \textbf{1036.61} & \textbf{1036.61} &    0.00 &  823.41 &   3 \\
            15 & 3 & 2 & 3 &          1014.04 & \textbf{1049.45} &    2.88 & 2054.35 &   2 & \textbf{1049.45} & \textbf{1049.45} &    0.00 & 1325.68 &   3 &          1021.59 & \textbf{1049.45} &    2.41 & 2419.87 &   1 \\
            15 & 3 & 3 & 2 &          1065.96 & \textbf{1114.68} &    4.39 &    tlim &   0 & \textbf{1114.68} & \textbf{1114.68} &    0.00 & 1104.10 &   3 &          1067.19 & \textbf{1114.68} &    4.28 & 2973.83 &   1 \\
            15 & 3 & 3 & 3 &          1101.81 & \textbf{1160.32} &    5.05 &    tlim &   0 & \textbf{1153.10} &          1160.74 &    0.65 & 1358.60 &   2 &          1096.70 & \textbf{1160.32} &    5.49 &    tlim &   0 \\
            15 & 4 & 2 & 2 & \textbf{1036.61} & \textbf{1036.61} &    0.00 & 1822.00 &   3 & \textbf{1036.61} & \textbf{1036.61} &    0.00 &  510.81 &   3 & \textbf{1036.61} & \textbf{1036.61} &    0.00 &  825.53 &   3 \\
            15 & 4 & 2 & 3 &          1081.47 & \textbf{1084.12} &    0.28 & 2439.70 &   2 & \textbf{1084.12} & \textbf{1084.12} &    0.00 & 1108.56 &   3 & \textbf{1084.12} & \textbf{1084.12} &    0.00 & 1637.68 &   3 \\
            15 & 4 & 3 & 2 &          1073.81 &          1150.81 &    6.81 &    tlim &   0 & \textbf{1142.24} & \textbf{1149.80} &    0.69 & 1420.98 &   2 &          1082.75 &          1150.81 &    6.03 &    tlim &   0 \\
            15 & 4 & 3 & 3 &          1111.81 &          1203.22 &    7.72 &    tlim &   0 & \textbf{1190.64} & \textbf{1202.40} &    1.03 & 1430.51 &   2 &          1113.29 &          1204.85 &    7.73 &    tlim &   0 \\
            \cmidrule(lr){1-4} \cmidrule(lr){5-9} \cmidrule(lr){10-14} \cmidrule(lr){15-19}
            \multicolumn{4}{l}{sum/avg (15)} & 1065.27 & 1104.48 & 3.39 & 2815.96 & 10 & \textbf{1100.93} & \textbf{1104.30} & 0.30 & 1096.12 & 21 & 1067.36 & 1104.68 & 3.24 & 2435.64 & 11 \\
            \cmidrule(lr){1-4} \cmidrule(lr){5-9} \cmidrule(lr){10-14} \cmidrule(lr){15-19}
            20 & 2 & 2 & 2 &          1225.52 & \textbf{1275.44} &    3.87 & 3274.98 &   1 & \textbf{1243.73} &          1277.82 &    2.61 & 1354.86 &   2 &          1228.48 & \textbf{1275.44} &    3.62 & 2266.04 &   2 \\
            20 & 2 & 2 & 3 &          1262.37 &          1330.87 &    5.14 &    tlim &   0 & \textbf{1294.95} & \textbf{1329.12} &    2.56 & 1458.53 &   2 &          1251.59 &          1332.91 &    6.09 &    tlim &   0 \\
            20 & 3 & 2 & 2 &          1197.44 & \textbf{1269.71} &    5.60 &    tlim &   0 & \textbf{1243.52} & \textbf{1269.71} &    2.02 & 1271.07 &   2 &          1236.21 & \textbf{1269.71} &    2.58 & 1531.30 &   2 \\
            20 & 3 & 2 & 3 &          1220.35 & \textbf{1301.73} &    6.13 &    tlim &   0 & \textbf{1301.73} & \textbf{1301.73} &    0.00 & 1127.63 &   3 &          1253.86 & \textbf{1301.73} &    3.58 & 2510.73 &   2 \\
            20 & 3 & 3 & 2 &          1210.01 & \textbf{1289.02} &    5.94 &    tlim &   0 & \textbf{1260.69} &          1296.61 &    2.60 & 2661.64 &   1 &          1215.65 &          1295.67 &    5.95 &    tlim &   0 \\
            20 & 3 & 3 & 3 &          1223.85 & \textbf{1321.84} &    7.26 &    tlim &   0 & \textbf{1281.30} & \textbf{1321.84} &    2.97 &    tlim &   0 &          1204.51 &          1337.03 &    9.69 &    tlim &   0 \\
            20 & 5 & 3 & 2 &          1200.10 &          1280.16 &    6.12 &    tlim &   0 & \textbf{1241.98} &          1280.43 &    2.91 & 2986.79 &   1 &          1235.59 & \textbf{1278.32} &    3.24 & 2739.77 &   1 \\
            20 & 5 & 3 & 3 &          1228.58 &          1324.16 &    7.12 &    tlim &   0 & \textbf{1270.11} & \textbf{1322.70} &    3.84 & 2863.57 &   1 &          1251.23 &          1402.59 &    9.43 & 2891.87 &   1 \\
            \cmidrule(lr){1-4} \cmidrule(lr){5-9} \cmidrule(lr){10-14} \cmidrule(lr){15-19}
            \multicolumn{4}{l}{sum/avg (20)} & 1221.03 & \textbf{1299.11} & 5.90 & 3559.37 & 1 & \textbf{1267.25} & 1300.00 & 2.44 & 2165.51 & 12 & 1234.64 & 1311.67 & 5.52 & 2842.59 & 8 \\
            \cmidrule(lr){1-4} \cmidrule(lr){5-9} \cmidrule(lr){10-14} \cmidrule(lr){15-19}
            \multicolumn{4}{l}{overall sum/avg} & 1041.02 & \textbf{1083.67} & 3.38 & 2333.13 & 29 & \textbf{1070.80} & 1083.93 & {0.99} & 1190.99 & 51 & 1046.73 & 1088.31 & 3.19 & 1942.39 & 37 \\
            \bottomrule
        \end{tabular}
    }
\end{sidewaystable}

From Table~\ref{tab:math_models_comparison}, we can observe that just on a few large-size instances the time-based model and the node-based model find better results in terms of average upper bound. Overall, the flow-based model outperforms the other two models in terms of average lower bound, average percentage gap, average run time, and number of optimal solutions obtained. Consequently, we adopted this model to assess the quality of the solutions obtained by the ILS (see Section \ref{subsec:ILS_parameter_tuning}).

For all the instances belonging to the medium- and large-size sets, the mathematical models could not obtain proven optimal solutions and the computer frequently ran out of memory because of the large model size. Overall, we can conclude that the results prove the need of a good heuristic for these instances. This need is further motivated by the dimension of the original real-world problem, where the number of visits per day (i.e., around 70) is out of scale if compared to the size of instances solved to optimality within the time limit.

To assess the performance of the proposed valid inequalities, six small-size instances were selected and solved running the three models with and without the addition of the valid inequalities. The results are reported in Table~\ref{tab:valid_inequalities}. We can notice that the inequalities help improve the performance of all models, by reducing the average percentage gap and execution time, and increasing the number of proven optimal solutions.

\begin{table}
    \caption{\small Effect of valid inequalities on six small-size instances} \label{tab:valid_inequalities}
    \centering
    \scriptsize
    \setlength{\tabcolsep}{3.5pt}
    \adjustbox{max width=\columnwidth}{
        \begin{tabular}{l r r r r r r r r r r}
            \toprule
             & \multicolumn{5}{c}{without valid inequalities} & \multicolumn{5}{c}{with valid inequalities} \\
            \multicolumn{1}{l}{mathematical model}  & $z_{lb}$ & $z_{ub}$ & $\%$gap &    t(s) & opt & $z_{lb}$ & $z_{ub}$ & $\%$gap &    t(s) & opt \\
            \cmidrule(lr){1-1} \cmidrule(lr){2-6} \cmidrule(lr){7-11}
            \multicolumn{1}{l}{time-based} & 1036.25 & 1077.49 & 3.24 & 2484.51 & 2 & 1040.67 & 1076.85 & 2.81 & 2288.63 & 3 \\
            \multicolumn{1}{l}{flow-based} & 1062.79 & 1076.85 & 1.00 &  953.69 & 5 & 1064.43 & 1076.85 & 0.88 & 751.28 & 5 \\
            \multicolumn{1}{l}{node-based} & 1033.41 & 1084.99 & 3.86 & 1912.63 & 3 & 1037.93 & 1068.09 & 2.38 & 1745.36 & 4 \\
            \bottomrule
        \end{tabular}
    }
\end{table}

\subsection{ILS Parameter Tuning} \label{subsec:ILS_parameter_tuning}
The ILS procedure adopts four {main} parameters (i.e., $\alpha$, $\beta$, $\gamma$ and $max_{iter}$). To set their values, we randomly selected six instances (two with $0 \leq n \leq 20$, two with $50 \leq n \leq 100$, and two with $n=200$). We then tested the ILS on these instances by attempting all possible combinations of parameter values chosen in the sets 
$\alpha \in \{0.05. 0.10, 0.15, 0.25\}$, 
$\beta \in \{2, 5, 10, 20\}$, 
$\gamma \in \{50, 100\}$ and
$max_{iter} \in \{200, 500, 1000, 5000\}$. 
The results are reported in Table~\ref{tab:param_tuning}. For each combination of parameters, column ``t(s)'' gives the average ILS run time on the six instances, and column ``$\%$gap'' gives the average gap computed as the average over the six instances of $100(z-z^*)/z^*$. Here, $z$ is the value of the solution obtained by the given configuration and $z^*$ is the value of the best solution obtained by all configurations.

{The configuration with $\alpha=0.10$, $\beta=5$, $\gamma=50$ and $max_{iter}=1000$ is the one that obtained the best results (highlighted in bold in the table). It could always achieve the best solution values, at the expense of a limited increase in the computing time with respect to configurations adopting a smaller number of iterations. This configuration was thus adopted for all successive ILS tests.}

\begin{table}
    \caption{\small ILS parameter tuning. Best configuration in \textbf{boldface}} \label{tab:param_tuning}
    \centering
    \scriptsize
    \setlength{\tabcolsep}{2.5pt}
    \adjustbox{max width=\columnwidth}{
    \begin{tabular}{l r r r r r r r r r r r r r r r r}
        \toprule
         & \multicolumn{16}{c}{$(\gamma, max_{iter})$} \\
        \cmidrule(lr){2-17}
         & \multicolumn{2}{c}{$(50, 200)$} & \multicolumn{2}{c}{$(50, 500)$} & \multicolumn{2}{c}{$(50, 1000)$} & \multicolumn{2}{c}{$(50, 5000)$} & \multicolumn{2}{c}{$(100, 200)$} & \multicolumn{2}{c}{$(100, 500)$} & \multicolumn{2}{c}{$(100, 1000)$} & \multicolumn{2}{c}{$(100, 5000)$} \\
        $(\alpha, \beta)$ & {t(s)} & {$\%$gap} & {t(s)} & {$\%$gap} & {t(s)} & {$\%$gap} & {t(s)} & {$\%$gap} & {t(s)} & {$\%$gap} & {t(s)} & {$\%$gap} & {t(s)} & {$\%$gap} & {t(s)} & {$\%$gap} \\
        \cmidrule(lr){1-1} \cmidrule(lr){2-3} \cmidrule(lr){4-5} \cmidrule(lr){6-7} \cmidrule(lr){8-9} \cmidrule(lr){10-11} \cmidrule(lr){12-13} \cmidrule(lr){14-15} \cmidrule(lr){16-17}
        (0.05,2)  & 1.53 & 0.91 & 1.79 & 0.87 &          2.13 &          0.84 & 2.79 & 0.82 & 1.56 & 0.83 & 1.80 & 0.79 & 1.89 & 0.78 & 3.28 & 0.77 \\
        (0.05,5)  & 1.67 & 0.76 & 1.83 & 0.75 &          2.27 &          0.73 & 2.91 & 0.72 & 1.79 & 0.75 & 1.90 & 0.74 & 2.02 & 0.74 & 3.49 & 0.73 \\
        (0.05,10) & 1.91 & 0.76 & 2.18 & 0.74 &          2.66 &          0.73 & 3.41 & 0.72 & 2.08 & 0.73 & 2.19 & 0.73 & 2.26 & 0.73 & 3.91 & 0.71 \\
        (0.05,20) & 1.73 & 0.76 & 1.93 & 0.74 &          2.34 &          0.73 & 2.86 & 0.72 & 2.20 & 0.73 & 2.35 & 0.72 & 2.43 & 0.69 & 4.67 & 0.69 \\
        \cmidrule(lr){1-1} \cmidrule(lr){2-3} \cmidrule(lr){4-5} \cmidrule(lr){6-7} \cmidrule(lr){8-9} \cmidrule(lr){10-11} \cmidrule(lr){12-13} \cmidrule(lr){14-15} \cmidrule(lr){16-17}
        (0.10,2)  & 2.09 & 0.08 & 2.33 & 0.03 &          2.68 &          0.02 & 3.58 & 0.01 & 1.91 & 0.13 & 1.97 & 0.11 & 2.09 & 0.10 & 2.58 & 0.10 \\
        (0.10,5)  & 2.74 & 0.08 & 3.25 & 0.02 & \textbf{4.02} & \textbf{0.00} & 5.44 & 0.00 & 2.06 & 0.11 & 2.38 & 0.07 & 2.89 & 0.06 & 3.28 & 0.05 \\
        (0.10,10) & 3.13 & 0.08 & 4.06 & 0.02 &          5.04 &          0.00 & 6.47 & 0.00 & 2.49 & 0.07 & 2.61 & 0.07 & 3.05 & 0.06 & 3.39 & 0.05 \\
        (0.10,20) & 3.57 & 0.08 & 4.53 & 0.02 &          5.23 &          0.00 & 8.02 & 0.00 & 2.84 & 0.07 & 3.11 & 0.06 & 3.24 & 0.06 & 3.46 & 0.05 \\
        \cmidrule(lr){1-1} \cmidrule(lr){2-3} \cmidrule(lr){4-5} \cmidrule(lr){6-7} \cmidrule(lr){8-9} \cmidrule(lr){10-11} \cmidrule(lr){12-13} \cmidrule(lr){14-15} \cmidrule(lr){16-17}
        (0.15,2)  & 1.83 & 0.43 & 2.06 & 0.39 &          2.30 &          0.38 & 3.12 & 0.35 & 2.04 & 0.38 & 2.37 & 0.35 & 2.49 & 0.35 & 3.12 & 0.35 \\
        (0.15,5)  & 2.49 & 0.40 & 2.86 & 0.38 &          3.13 &          0.35 & 5.08 & 0.35 & 2.33 & 0.36 & 2.59 & 0.35 & 3.20 & 0.34 & 4.85 & 0.33 \\
        (0.15,10) & 3.35 & 0.40 & 3.88 & 0.38 &          4.16 &          0.35 & 6.37 & 0.35 & 2.48 & 0.35 & 3.79 & 0.33 & 4.11 & 0.33 & 5.09 & 0.33 \\
        (0.15,20) & 4.55 & 0.40 & 5.02 & 0.38 &          5.23 &          0.35 & 8.64 & 0.35 & 2.71 & 0.35 & 4.26 & 0.33 & 4.82 & 0.33 & 5.94 & 0.32 \\
        \cmidrule(lr){1-1} \cmidrule(lr){2-3} \cmidrule(lr){4-5} \cmidrule(lr){6-7} \cmidrule(lr){8-9} \cmidrule(lr){10-11} \cmidrule(lr){12-13} \cmidrule(lr){14-15} \cmidrule(lr){16-17}
        (0.25,2)  & 2.25 & 1.34 & 2.64 & 1.07 &          3.30 &          0.94 & 4.56 & 0.92 & 2.21 & 0.88 & 2.27 & 0.86 & 2.84 & 0.86 & 4.19 & 0.86 \\
        (0.25,5)  & 2.54 & 1.18 & 2.93 & 0.91 &          4.05 &          0.89 & 5.62 & 0.89 & 2.68 & 0.86 & 3.16 & 0.85 & 3.74 & 0.85 & 4.80 & 0.83 \\
        (0.25,10) & 3.00 & 1.16 & 3.94 & 0.90 &          4.94 &          0.86 & 7.33 & 0.86 & 3.52 & 0.83 & 4.86 & 0.82 & 6.07 & 0.82 & 7.83 & 0.82 \\
        (0.25,20) & 3.72 & 1.16 & 5.12 & 0.90 &          6.31 &          0.86 & 8.89 & 0.86 & 4.67 & 0.83 & 6.13 & 0.82 & 6.63 & 0.82 & 8.46 & 0.82 \\
        \bottomrule
    \end{tabular}
    }
\end{table}

\subsection{ILS Evaluation} \label{subsec:ILS_evaluation}
In this section, we investigate the performance of the ILS. In Table~\ref{tab:comparison_best_model_ILS_small}, the results of the ILS are compared with those obtained by the best mathematical model (i.e., the flow-based one) on groups of three instances per line. We recall that column ``$z_{ub}$'' gives the average upper bound value, column ``opt'' the  number of proven optimal solutions, and column ``t(s)'' the average run time. The ILS was executed five times on each instance. We report the best, average and worst solution values achieved, as well as their standard deviation, in columns ``$z_{best}$'', ``$z_{avg}$'', ``$z_{worst}$'' and ``$\sigma_z$'', respectively. More in detail, $z_{best}$ gives the average of the best solution values produced on the three instances, $z_{avg}$ the average of the average values, and $z_{worst}$ the average of the worst values. The average computational time is shown in column ``t(s)''.

According to the results, for those groups of three instances that were all solved to optimality by the flow-based model, the ILS obtained the same optimal values in a shorter computational time. For all the remaining small-size sets, the ILS achieved better values than the flow-based model (without proof of their optimality). In addition, the constantly null average standard deviation among the different runs indicates the robustness of the algorithm on these very simple instances. When comparing the average run times, we can notice that the ILS needed an overall average time of just 0.23 seconds against the 1,190.99 seconds of the flow-based model.

\begin{table}
  \caption{\small Computational results on small-size instances (three inst. per line)} \label{tab:comparison_best_model_ILS_small}
  \centering
  \scriptsize
  \adjustbox{max width=\columnwidth}{
    \begin{tabular}{c c c c r r r r r r r r}
        \toprule
           &   &   &   & \multicolumn{3}{c}{flow-based} & \multicolumn{5}{c}{ILS} \\
        \cmidrule(lr){5-7} \cmidrule(lr){8-12}
        $|V_1 \cup V_2|$ & $|V_2|$ & $|V_3|$ & $|K|$ & $z_{ub}$ & t(s) & opt & $z_{best}$ & $z_{avg}$ & $z_{worst}$ & $\sigma_z$ & t(s) \\
        \cmidrule(lr){1-4} \cmidrule(lr){5-7} \cmidrule(lr){8-12}
        10 & 1 & 1 & 1 &  706.39 &    1.13 & 3 &  706.39 &  706.39 &   706.39 &     0.00 & 0.00 \\
        10 & 1 & 1 & 2 &  730.91 &    3.26 & 3 &  730.91 &  730.91 &   730.91 &     0.00 & 0.00 \\
        10 & 2 & 1 & 1 &  743.25 &    3.13 & 3 &  743.25 &  743.25 &   743.25 &     0.00 & 0.00 \\
        10 & 2 & 1 & 2 &  793.75 &    9.45 & 3 &  793.75 &  793.75 &   793.75 &     0.00 & 0.00 \\
        10 & 2 & 2 & 1 &  803.96 &   16.73 & 3 &  803.96 &  803.96 &   803.96 &     0.00 & 0.00 \\
        10 & 2 & 2 & 2 &  833.78 &   75.01 & 3 &  833.78 &  833.78 &   833.78 &     0.00 & 0.00 \\
        \cmidrule(lr){1-4} \cmidrule(lr){5-7} \cmidrule(lr){8-12}
        \multicolumn{4}{l}{sum/avg (10)} & 768.67 & 18.12 & 18 & 768.67 & 768.67 & 768.67 & 0.00 & 0.00 \\
        \cmidrule(lr){1-4} \cmidrule(lr){5-7} \cmidrule(lr){8-12}
        15 & 3 & 2 & 2 & 1036.61 &  509.71 & 3 & 1036.61 & 1036.61 &  1036.61 &     0.00 & 0.14 \\
        15 & 3 & 2 & 3 & 1049.45 & 1325.68 & 3 & 1049.45 & 1049.45 &  1049.45 &     0.00 & 0.22 \\
        15 & 3 & 3 & 2 & 1114.68 & 1104.10 & 3 & 1114.68 & 1114.68 &  1114.68 &     0.00 & 0.14 \\
        15 & 3 & 3 & 3 & 1160.74 & 1358.60 & 2 & 1155.96 & 1155.96 &  1155.96 &     0.00 & 0.24 \\
        15 & 4 & 2 & 2 & 1036.61 &  510.81 & 3 & 1036.61 & 1036.61 &  1036.61 &     0.00 & 0.16 \\
        15 & 4 & 2 & 3 & 1084.12 & 1108.56 & 3 & 1084.12 & 1084.12 &  1084.12 &     0.00 & 0.22 \\
        15 & 4 & 3 & 2 & 1149.80 & 1420.98 & 2 & 1146.56 & 1146.56 &  1146.56 &     0.00 & 0.20 \\
        15 & 4 & 3 & 3 & 1202.40 & 1430.51 & 2 & 1202.40 & 1202.40 &  1202.40 &     0.00 & 0.27 \\
        \cmidrule(lr){1-4} \cmidrule(lr){5-7} \cmidrule(lr){8-12}
        \multicolumn{4}{l}{sum/avg (15)} & 1104.30 & 1096.12 & 21 & 1103.30 & 1103.30 & 1103.30 & 0.00 & 0.20 \\
        \cmidrule(lr){1-4} \cmidrule(lr){5-7} \cmidrule(lr){8-12}
        20 & 2 & 2 & 2 & 1277.82 & 1354.86 & 2 & 1275.44 & 1275.44 &  1275.44 &     0.00 & 0.26 \\
        20 & 2 & 2 & 3 & 1329.12 & 1458.53 & 2 & 1316.03 & 1316.03 &  1316.03 &     0.00 & 0.33 \\
        20 & 3 & 2 & 2 & 1269.71 & 1271.07 & 2 & 1262.52 & 1262.52 &  1262.52 &     0.00 & 0.35 \\
        20 & 3 & 2 & 3 & 1301.73 & 1127.63 & 3 & 1301.73 & 1301.73 &  1301.73 &     0.00 & 0.41 \\
        20 & 3 & 3 & 2 & 1296.61 & 2661.64 & 1 & 1277.25 & 1277.25 &  1277.25 &     0.00 & 0.35 \\
        20 & 3 & 3 & 3 & 1321.84 &    tlim & 0 & 1301.37 & 1301.37 &  1301.37 &     0.00 & 0.61 \\
        20 & 5 & 3 & 2 & 1280.43 & 2986.79 & 1 & 1265.46 & 1265.46 &  1265.46 &     0.00 & 0.33 \\
        20 & 5 & 3 & 3 & 1322.70 & 2863.57 & 1 & 1303.93 & 1303.93 &  1303.93 &     0.00 & 0.73 \\
        \cmidrule(lr){1-4} \cmidrule(lr){5-7} \cmidrule(lr){8-12}
        \multicolumn{4}{l}{sum/avg (20)} & 1300.00 & 2165.51 & 12 & 1287.97 & 1287.97 & 1287.97 & 0.00 & 0.42 \\
        \cmidrule(lr){1-4} \cmidrule(lr){5-7} \cmidrule(lr){8-12}
        \multicolumn{4}{l}{overall sum/avg} & 1083.93 & 1190.99 & 51 & 1079.19 & 1079.19 & 1079.19 & 0.00 & 0.23 \\
        \bottomrule
    \end{tabular}
  }
\end{table}

In Table~\ref{tab:results_ILS_medium_large}, we report the results of the ILS on medium- and large-size instances. On instances having $|V_1 \cup V_2| = 50$ the average standard deviation is 0.00, on those having $|V_1 \cup V_2| = 100$ it becomes 0.50, while on those having $|V_1 \cup V_2| = 200$ it increases to 0.92, thus resulting in an overall average standard deviation of 0.47. This confirms the robustness of the algorithm. Concerning the run time, the ILS took on average 1.91 seconds to solve instances having $|V_1 \cup V_2| = 50$, 8.60 seconds for those having $|V_1 \cup V_2| = 100$, and 13.02 seconds for those having $|V_1 \cup V_2| = 200$. The overall average run time is 7.84 seconds, proving that the method is suitable for a quick use in practical situations.

\begin{table}
  \caption{\small Computational results on medium- and large-size instances (three inst. per line)} \label{tab:results_ILS_medium_large}
  \centering
  \scriptsize
  \setlength{\tabcolsep}{11pt}
  \adjustbox{max width=\columnwidth}{
    \begin{tabular}{c c c c r r r r r}
        \toprule
           &   &   &   & \multicolumn{5}{c}{ILS} \\
        \cmidrule(lr){5-9}
        $|V_1 \cup V_2|$ & $|V_2|$ & $|V_3|$ & $|K|$ & $z_{best}$ & $z_{avg}$ & $z_{worst}$ & $\sigma_z$ & t(s) \\
        \cmidrule(lr){1-4} \cmidrule(lr){5-9}
         50 &  5 &  5 &  5 & 2583.27 & 2583.27 & 2583.27 & 0.00 &  1.17 \\
         50 &  5 &  5 &  8 & 2721.90 & 2721.90 & 2721.90 & 0.00 &  1.59 \\
         50 &  8 &  8 &  5 & 2883.07 & 2883.07 & 2883.07 & 0.00 &  1.67 \\
         50 &  8 &  8 &  8 & 3001.70 & 3001.70 & 3001.70 & 0.00 &  2.10 \\
         50 & 10 &  5 &  5 & 2664.02 & 2664.02 & 2664.02 & 0.00 &  2.09 \\
         50 & 10 &  5 &  8 & 2807.41 & 2807.41 & 2807.41 & 0.00 &  2.19 \\
         50 & 10 &  8 &  5 & 2863.64 & 2863.64 & 2863.64 & 0.00 &  1.91 \\
         50 & 10 &  8 &  8 & 3003.07 & 3003.07 & 3003.07 & 0.00 &  2.52 \\
        \cmidrule(lr){1-4} \cmidrule(lr){5-9}
        \multicolumn{4}{l}{avg (50)} & 2816.01 & 2816.01 & 2816.01 & 0.00 & 1.91 \\
        \cmidrule(lr){1-4} \cmidrule(lr){5-9}
        100 &  5 &  5 & 10 & 4430.65 & 4430.92 & 4431.53 & 0.40 &  7.61 \\
        100 &  5 &  5 & 15 & 4642.24 & 4642.45 & 4643.13 & 0.39 &  8.29 \\
        100 & 10 &  5 & 10 & 4507.07 & 4507.27 & 4508.07 & 0.45 &  7.19 \\
        100 & 10 &  5 & 15 & 4750.73 & 4750.92 & 4751.65 & 0.41 &  8.17 \\
        100 & 10 & 10 & 10 & 4856.94 & 4857.16 & 4857.99 & 0.47 &  9.03 \\
        100 & 10 & 10 & 15 & 5062.41 & 5062.62 & 5063.43 & 0.45 &  9.43 \\
        100 & 15 & 10 & 10 & 4826.28 & 4826.62 & 4827.96 & 0.75 &  9.18 \\
        100 & 15 & 10 & 15 & 5070.19 & 5070.50 & 5071.74 & 0.69 &  9.90 \\
        \cmidrule(lr){1-4} \cmidrule(lr){5-9}
        \multicolumn{4}{l}{avg (100)} & 4768.31 & 4768.56 & 4769.44 & 0.50 & 8.60 \\
        \cmidrule(lr){1-4} \cmidrule(lr){5-9}
        200 & 10 & 10 & 15 & 8244.39 & 8244.97 & 8246.12 & 0.82 &  9.86 \\
        200 & 10 & 10 & 20 & 8636.53 & 8637.11 & 8638.33 & 0.84 & 10.34 \\
        200 & 20 & 10 & 15 & 8550.63 & 8551.32 & 8552.62 & 0.98 & 12.27 \\
        200 & 20 & 10 & 20 & 8814.41 & 8815.00 & 8816.05 & 0.82 & 13.29 \\
        200 & 20 & 20 & 15 & 9128.90 & 9129.63 & 9130.63 & 0.82 & 13.66 \\
        200 & 20 & 20 & 20 & 9305.35 & 9305.98 & 9307.13 & 0.81 & 14.96 \\
        200 & 30 & 20 & 15 & 9372.60 & 9373.86 & 9375.17 & 1.16 & 14.05 \\
        200 & 30 & 20 & 20 & 9497.20 & 9498.20 & 9499.67 & 1.10 & 15.70 \\
        \cmidrule(lr){1-4} \cmidrule(lr){5-9}
        \multicolumn{4}{l}{avg (200)} & 8943.75 & 8944.51 & 8945.72 & 0.92 & 13.02 \\
        \cmidrule(lr){1-4} \cmidrule(lr){5-9}
        \multicolumn{4}{l}{overall avg} & 5509.36 & 5509.69 & 5510.39 & 0.47 & 7.84 \\
        \bottomrule
    \end{tabular}
  }
\end{table}

\begin{table}
	\caption{\small Percentage of the computational time needed by each ILS component} \label{tab:ILS_component}
    \centering
    \scriptsize
    \setlength{\tabcolsep}{10pt}
    \adjustbox{max width=\columnwidth}{
        \begin{tabular}{l r r r r r r r}
    	    \toprule
    	    Set         &     LS1 &     LS2 &     LS3 &     LS4 &     LS5 &     S1 &     S2 \\
    	    \cmidrule(lr){1-1} \cmidrule(lr){2-8}
    	    Small-size  &  3.63\% & 39.41\% & 28.65\% &  6.37\% & 20.50\% & 0.37\% & 1.07\% \\
    	    Medium-size & 59.50\% &  9.02\% &  1.30\% & 18.36\% & 10.69\% & 0.23\% & 0.89\% \\
    	    Large-size  & 50.03\% &  3.48\% &  3.33\% & 32.35\% &  9.98\% & 0.13\% & 0.71\% \\ 
    	    \bottomrule
        \end{tabular}
    }
\end{table}

Finally, Table~\ref{tab:ILS_component} reports a sensitivity analysis on the average percentage of computational time needed by each ILS component, grouped by set of instances. On the small-size sets, LS2 and LS3 are the most time-consuming local search procedures, while for medium- and large-size sets the largest effort is required by LS1 and LS4.

\subsection{Results on Realistic Instances} \label{subsec:realistic_instances}
{The flow-based model and the ILS were also tested on a set of realistic instances generated from the WDN in the city of Mashhad (Iran). Our real case study consists of 3,124 households/shops, 293 reservoirs/tanks, 356 wells {and 14 treatment plants}. For all of these nodes the exact locations were collected.}

Following the same rationale described in Section~\ref{subsec:test_instances}, we generated 108 realistic instances divided into two sets of \emph{small-size} and \emph{medium- and large-size} instances, each comprising different subsets having homogeneous values of $|V_1 \cup V_2|$, $(|V_2|,|V_3|)$, and $|K|$, and composed by three random instances per subset. The resulting sets are:
\begin{itemize}
\item \emph{Small-size}: 
12 instances with $|V_1 \cup V_2|$=10, 
$(|V_2|,|V_3|) \in \{(1,1)$, $(2,1)$, $(2,2)\}$, and $|K| \in \{1,2\}$;
12 instances with $|V_1 \cup V_2|$=15,
$(|V_2|,|V_3|) \in \{(1,1)$, $(2,1)$, $(2,2)\}$, and $|K| \in \{1, 2, 3\}$;
12 instances with $|V_1 \cup V_2|$=20,
$(|V_2|,|V_3|) \in \{(1,1)$, $(2,1)$, $(2,2)\}$, and $|K| \in \{2, 3\}$;
\item \emph{Medium- and large-size}:
12 instances with $|V_1 \cup V_2|$=40,
$(|V_2|,|V_3|) \in \{(4,2)$, $(4,3)$, $(6,2)$, $(6,3)\}$, and $|K| \in \{2, 3\}$;
12 instances with $|V_1 \cup V_2|$=50,
$(|V_2|,|V_3|)$ $\in$ $\{(4,2)$, $(4,3)$, $(6,2)$, $(6,3)\}$, and $|K| \in \{2, 3\}$;
12 instances with $|V_1 \cup V_2|$=60,
$(|V_2|,|V_3|) \in \{(4,2)$, $(4,3)$, $(6,2)$, $(6,3)\}$, and $|K| \in \{2, 3\}$;
12 instances with $|V_1 \cup V_2|$=100,
$(|V_2|,|V_3|) \in \{(8,4)$, $(8,5)$, $(10,4)$, $(10,5)\}$, and $|K| \in \{4, 5\}$;
12 instances with $|V_1 \cup V_2|$=150,
$(|V_2|,|V_3|) \in \{(8,4)$, $(8,5)$, $(10,4)$, $(10,5)\}$, and $|K| \in \{4, 5\}$;
12 instances with $|V_1 \cup V_2|$=200,
$(|V_2|,|V_3|) \in \{(8,4)$, $(8,5)$, $(10,4)$, $(10,5)\}$, and $|K| \in \{4, 5\}$.
\end{itemize}

{For each instance, the coordinates of {the} nodes were randomly selected among the given real locations.
The flow-based model and the ILS were used to run the experiments. The results are reported in the Tables~\ref{tab:comparison_best_model_ILS_realistic_small} and \ref{tab:results_ILS_realistic_medium_large}. In Table~\ref{tab:comparison_best_model_ILS_realistic_small}, the results of the ILS are compared with those obtained by the flow-based model. We recall that columns ``$z_{ub}$'', ``t(s)'' and ``opt'' give the average upper bound value, the average run time and the total number of instances solved to proven optimality by the mathematical model, respectively. Note that an entry ``tlim'' indicates that the time limit of 3,600 CPU seconds was reached for all the three instances in the group. Conversely, columns ``$z_{best}$'', ``$z_{avg}$'', ``$z_{worst}$'', ``$\sigma_z$'' and ``t(s)'' give the best, average and worst solution values, the standard deviation and the computational time of the ILS, respectively.}

We can notice that on small-size instances, the flow-based model and the ILS obtained the same optimal values on instances having $|V_1 \cup V_2| \in \{10, 15\}$, and on one subset out of four of instances having $|V_1 \cup V_2| = 20$. For the remaining subsets, the ILS achieved better values than the flow-based model (again, without proof of their optimality).

\begin{table}
    \caption{\small Computational results on realistic small-size instances (three inst. per line)}
    \label{tab:comparison_best_model_ILS_realistic_small}
    \centering
    \scriptsize
    \setlength{\tabcolsep}{5pt}
    \adjustbox{max width=\columnwidth}{
        \begin{tabular}{c c c c r r r r r r r r}
            \toprule
               &   &   &   & \multicolumn{3}{c}{flow-based} & \multicolumn{5}{c}{ILS} \\
            \cmidrule(lr){5-7} \cmidrule(lr){8-12}
            $|V_1 \cup V_2|$ & $|V_2|$ & $|V_3|$ & $|K|$ & $z_{ub}$ & t(s) & opt & $z_{best}$ & $z_{avg}$ & $z_{worst}$ & $\sigma_z$ & t(s) \\
            \cmidrule(lr){1-4} \cmidrule(lr){5-7} \cmidrule(lr){8-12}
            10 & 1 & 1 & 1 &  82475.75 &   12.69 & 3 &  82475.75 &  82475.75 &  82475.75 & 0.00 & 0.00 \\
            10 & 1 & 1 & 2 &  80629.79 &    8.95 & 3 &  80629.79 &  80629.79 &  80629.79 & 0.00 & 0.00 \\
            10 & 2 & 1 & 1 & 121874.44 &    5.53 & 3 & 121874.44 & 121874.44 & 121874.44 & 0.00 & 0.00 \\
            10 & 2 & 2 & 2 & 116649.98 &  135.45 & 3 & 116649.98 & 116649.98 & 116649.98 & 0.00 & 0.00 \\
            \cmidrule(lr){1-4} \cmidrule(lr){5-7} \cmidrule(lr){8-12}
            \multicolumn{4}{l}{sum/avg (10)} & 100407.49 & 40.66 & 12 & 100407.49 & 100407.49 & 100407.49 & 0.00 & 0.00 \\
            \cmidrule(lr){1-4} \cmidrule(lr){5-7} \cmidrule(lr){8-12}
            15 & 1 & 1 & 1 &  97874.84 & 1218.80 & 2 &  97874.84 &  97874.84 &  97874.84 & 0.00 & 0.05 \\
            15 & 1 & 1 & 2 & 186943.98 &   46.03 & 3 & 186943.98 & 186943.98 & 186943.98 & 0.00 & 0.08 \\
            15 & 2 & 1 & 2 & 126742.54 & 1204.46 & 2 & 126742.54 & 126742.54 & 126742.54 & 0.00 & 0.12 \\
            15 & 2 & 2 & 3 & 135254.03 & 2437.79 & 1 & 135254.03 & 135254.03 & 135254.03 & 0.00 & 0.12 \\
            \cmidrule(lr){1-4} \cmidrule(lr){5-7} \cmidrule(lr){8-12}
            \multicolumn{4}{l}{sum/avg (15)} & 136703.85 & 1226.77 & 8 & 136703.85 & 136703.85 & 136703.85 & 0.00 & 0.09 \\
            \cmidrule(lr){1-4} \cmidrule(lr){5-7} \cmidrule(lr){8-12}
            20 & 1 & 1 & 2 & 158682.88 & 1205.87 & 2 & 158588.73 & 158588.73 & 158588.73 & 0.00 & 0.18 \\
            20 & 1 & 1 & 3 & 175830.70 & 3342.99 & 1 & 175655.43 & 175655.43 & 175655.43 & 0.00 & 0.20 \\
            20 & 2 & 1 & 2 & 170015.90 & 1466.86 & 3 & 170015.90 & 170015.90 & 170015.90 & 0.00 & 0.20 \\
            20 & 2 & 2 & 3 & 188935.98 &    tlim & 0 & 187949.09 & 187949.09 & 187949.09 & 0.00 & 0.26 \\
            \cmidrule(lr){1-4} \cmidrule(lr){5-7} \cmidrule(lr){8-12}
            \multicolumn{4}{l}{sum/avg (20)} & 173366.36 & 2403.93 & 6 & 173052.29 & 173052.29 & 173052.29 & 0.00 & 0.21 \\
            \cmidrule(lr){1-4} \cmidrule(lr){5-7} \cmidrule(lr){8-12}
            \multicolumn{4}{l}{overall sum/avg} & 136825.90 & 1223.79 & 26 & 136721.21 & 136721.21 & 136721.21 & 0.00 & 0.10 \\
            \bottomrule
        \end{tabular}
    }
\end{table}

On medium- and large-size instances, the ILS  achieved very robust results on instances having $|V_1 \cup V_2| \in \{40, 50, 60\}$. Indeed, the standard deviation is constantly null for all the subgroups, and the run times are very short. 
The robustness of the ILS slightly decreases for instances having $|V_1 \cup V_2| \in \{100, 150, 200\}$, however remaining acceptable for a practical use. For these instances, the average run times are around 7.32, 9.80 and 11.08 seconds, respectively, thus confirming that the algorithm could efficiently solve realistic instances having a considerable number of nodes in a few seconds.

\begin{table}
    \caption{\small Comp. results on realistic medium- and large-size instances (three inst. per line)} \label{tab:results_ILS_realistic_medium_large}
    \centering
    \scriptsize
    \setlength{\tabcolsep}{9.5pt}
    \adjustbox{max width=\columnwidth}{
        \begin{tabular}{c c c c r r r r r}
            \toprule
               &   &   &   & \multicolumn{5}{c}{ILS} \\
            \cmidrule(lr){5-9}
            $|V_1 \cup V_2|$ & $|V_2|$ & $|V_3|$ & $|K|$ & $z_{best}$ & $z_{avg}$ & $z_{worst}$ & $\sigma_z$ & t(s) \\
            \cmidrule(lr){1-4} \cmidrule(lr){5-9}
             40 &  4 & 2 & 2 & 180613.80 & 180613.80 & 180613.80 & 0.00 &  0.78 \\
             40 &  4 & 3 & 3 & 219227.54 & 219227.54 & 219227.54 & 0.00 &  0.79 \\
             40 &  6 & 2 & 2 & 195113.71 & 195113.71 & 195113.71 & 0.00 &  0.85 \\
             40 &  6 & 3 & 3 & 201338.25 & 201338.25 & 201338.25 & 0.00 &  0.84 \\
            \cmidrule(lr){1-4} \cmidrule(lr){5-9}
            \multicolumn{4}{l}{avg (40)} & 199073.33 & 199073.33 & 199073.33 & 0.00 & 0.82 \\
            \cmidrule(lr){1-4} \cmidrule(lr){5-9}
             50 &  4 & 2 & 2 & 250479.23 & 250479.23 & 250479.23 & 0.00 &  0.97 \\
             50 &  4 & 3 & 3 & 278476.01 & 278476.01 & 278476.01 & 0.00 &  0.97 \\
             50 &  6 & 2 & 2 & 263769.77 & 263769.77 & 263769.77 & 0.00 &  1.18 \\
             50 &  6 & 3 & 3 & 293179.73 & 293179.73 & 293179.73 & 0.00 &  1.28 \\
            \cmidrule(lr){1-4} \cmidrule(lr){5-9}
            \multicolumn{4}{l}{avg (50)} & 271476.19 & 271476.19 & 271476.19 & 0.00 & 1.10 \\
            \cmidrule(lr){1-4} \cmidrule(lr){5-9}
             60 &  4 & 2 & 2 & 263976.92 & 263976.92 & 263976.92 & 0.00 &  1.67 \\
             60 &  4 & 3 & 3 & 264029.28 & 264029.28 & 264029.28 & 0.00 &  1.73 \\
             60 &  6 & 2 & 2 & 222473.25 & 222473.25 & 222473.25 & 0.00 &  2.07 \\
             60 &  6 & 3 & 3 & 291979.18 & 291979.18 & 291979.18 & 0.00 &  2.35 \\
            \cmidrule(lr){1-4} \cmidrule(lr){5-9}
            \multicolumn{4}{l}{avg (60)} & 260614.66 & 260614.66 & 260614.66 & 0.00 & 1.96 \\
            \cmidrule(lr){1-4} \cmidrule(lr){5-9}
            100 &  8 & 4 & 4 & 399442.47 & 399442.63 & 399443.14 & 0.30 &  6.41 \\
            100 &  8 & 5 & 5 & 438923.85 & 438923.90 & 438924.05 & 0.09 &  7.38 \\
            100 & 10 & 4 & 4 & 344216.84 & 344217.11 & 344217.64 & 0.38 &  7.72 \\
            100 & 10 & 5 & 5 & 376473.23 & 376473.42 & 376473.86 & 0.28 &  7.78 \\
            \cmidrule(lr){1-4} \cmidrule(lr){5-9}
            \multicolumn{4}{l}{avg (100)} & 389764.10 & 389764.27 & 389764.67 & 0.26 & 7.32 \\
            \cmidrule(lr){1-4} \cmidrule(lr){5-9}
            150 &  8 & 4 & 4 & 438156.93 & 438157.42 & 438158.64 & 0.72 &  9.64 \\
            150 &  8 & 5 & 5 & 440416.07 & 440416.44 & 440417.48 & 0.61 &  8.85 \\
            150 & 10 & 4 & 4 & 466864.57 & 466864.99 & 466866.19 & 0.70 & 10.54 \\
            150 & 10 & 5 & 5 & 569988.53 & 569988.92 & 569990.11 & 0.68 & 10.17 \\
            \cmidrule(lr){1-4} \cmidrule(lr){5-9}
            \multicolumn{4}{l}{avg (150)} & 478856.53 & 478856.94 & 478858.11 & 0.68 & 9.80 \\
            \cmidrule(lr){1-4} \cmidrule(lr){5-9}
            200 &  8 & 4 & 4 & 491986.12 & 491986.81 & 491988.10 & 0.86 & 10.76 \\
            200 &  8 & 5 & 5 & 495553.12 & 495553.59 & 495555.14 & 0.88 &  9.82 \\
            200 & 10 & 4 & 4 & 646227.25 & 646227.71 & 646228.81 & 0.71 & 11.93 \\
            200 & 10 & 5 & 5 & 696373.63 & 696374.24 & 696375.29 & 0.73 & 11.80 \\
            \cmidrule(lr){1-4} \cmidrule(lr){5-9}
            \multicolumn{4}{l}{avg (200)} & 582535.03 & 582535.58 & 582536.83 & 0.79 & 11.08 \\
            \cmidrule(lr){1-4} \cmidrule(lr){5-9}
            \multicolumn{4}{l}{overall avg} & 363719.97 & 363720.16 & 363720.63 & 0.29 & 5.35 \\
            \bottomrule
        \end{tabular}
    }
\end{table}

\section{Conclusions} \label{sec:conclusions}
In this paper, we introduced a generalization of the well-{known} Vehicle Routing Problem (VRP), called VRP for Water Distribution Networks (VRPWDN), that includes precedence constraints among nodes and multiple visits to some of the nodes. The problem is NP-hard in the strong sense and, to the best of our knowledge, has not yet been applied in the context of distribution networks where regular inspections have to be performed to detect potential sources of contamination. To solve the VRPWDN, three alternative mathematical models (time-based, flow-based and node-based) were proposed, and an Iterated Local Search (ILS) algorithm was developed.

Extensive computational tests on randomly generated small-size instances were performed to compare the performance of the three mathematical models, showing that the flow-based model outperforms the other two in terms of solution quality and speed. On the same instances, the accuracy of the ILS in finding good-quality solutions in a short time was proved. The ILS was also used to perform a series of tests on randomly generated medium- and large-size instances with up to 200 nodes, confirming its efficacy and robustness.

Additional computational tests were executed on small-, medium-, and large-size realistic instances derived from the Mashhad (Iran) distribution network, proving that our methods can be applied with profit even in a practical case.

Interesting future research directions include the application of the developed techniques to other related VRPs with precedence constraints and multiple visits. In addition, we are interested in studying the generalization of the VRPWDN to the case of multiple periods. In this generalization, one should first of all determine in which day inspecting the given nodes, and then creating the routes for each day.

\section*{Acknowledgements}
This work was supported in part by Research Deputy of Ferdowsi University of Mashhad, under Grant No. 48185, and by University of Modena and Reggio Emilia, under grant FAR 2021. This support is gratefully acknowledged.

\bibliographystyle{apalike}
\bibliography{ref}

\end{document}